\documentclass[11pt]{article}
\usepackage{amsmath}
\usepackage{amssymb}
\usepackage{graphics}
\usepackage{graphicx}
\usepackage{array}
\newtheorem{theo}{\hspace*{\parindent}Theorem}

\newtheorem{lemma}{\hspace*{\parindent}Lemma}

\def\arctanh{\mathrm{arctanh}}

\def\O{{\cal O}}

\newcounter{theremark}
\newcommand{\rem}{\par\refstepcounter{theremark}\textbf{Remark \arabic{theremark}.} }

\title{Asymptotic approximations for the first incomplete elliptic integral near logarithmic singularity}
\author{D. Karp\footnote{Institute of Applied Mathematics, Vladivostok, Russia,
e-mail:\emph{dmkrp@yandex.ru}}~
 and S.\,M.\,Sitnik\footnote{Voronezh Institute of the Ministry of
Internal Affairs of the Russian Federation,
e-mail:\emph{mathsms@yandex.ru}}}
\date{}
\begin{document}
\maketitle
\begin{center}
\parbox{12cm}{
\small\textbf{Abstract.}  We find two convergent series expansions
for Legendre's first incomplete elliptic integral $F(\lambda,k)$
in terms of recursively computed elementary functions. Both
expansions are valid at every point of the unit square
$0<\lambda,k<1$. Truncated expansions yield asymptotic
approximations for $F(\lambda,k)$ as $\lambda$ and/or $k$ tend to
unity, including the case when logarithmic singularity
$\lambda=k=1$ is approached from any direction. Explicit error
bounds are given at every order of approximation. For the reader's
convenience we present explicit expressions for low-order
approximations and numerical examples to illustrate their
accuracy. Our derivation is based on rearrangements of some known
double series expansions, hypergeometric summation algorithms and
inequalities for hypergeometric functions.}
\end{center}

\bigskip

Keywords: \emph{Incomplete elliptic integral, series expansion,
asymptotic approximation, hypergeometric inequality}

\bigskip

MSC2000: 33E05, 33C75, 33F05.

\bigskip

\paragraph{1. Introduction.}

Legendre's incomplete elliptic integral (EI) of the first kind is
defined by \cite[(12.2.7)]{AS}:
\begin{equation}\label{eq:F-defined}
F(\lambda,k)=\int\limits_{0}^{\lambda}\frac{dt}{\sqrt{(1-t^2)(1-k^2t^2)}}.
\end{equation}
It is one of the three canonical forms given by Legendre in terms
of which all elliptic integrals can be expressed. We will only
consider the most important case $0\leq{k}\leq{1}$,
$0\leq{\lambda}\leq{1}$.

The subject of series expansions and asymptotic approximations for
the first incomplete elliptic integral has a long history. An
expansion given by E.L.\,Kaplan in 1948 \cite{Kaplan} implicitly
contained an asymptotic approximation for $F(\lambda,k)$ near the
singular point $\lambda=k=1$  (see (\ref{eq:CG1}) below). Soon
after Kaplan's paper various series expansions for the first
incomplete EI were given by B.\,Radon (1950) in \cite{Radon} and
R.P.\,Kelisky (1959) in \cite{Kelisky}. A bit later B.C.\,Carlson
showed in \cite{Carlson1} that $F(\lambda,k)$ can be expressed in
terms of Appell's first hypergeometric series $F_1$ (see
\cite{Bat1}), which automatically lead to several series
expansions through known transformation formulas for $F_1$. In the
same paper he noted that one can derive rapidly convergent
expansions by first expressing Legendre's incomplete EIs in a
different form. This form had later become known as symmetric
standard EIs. B.C.\,Carlson proved that instead of three
Legendre's EIs one can use three symmetric standard elliptic
integrals as canonical forms. The first symmetric standard
elliptic integral is defined by
\cite{Carlson1,CarlsonBook,CG1,CG2}:
\begin{equation}\label{eq:Rf}
R_{F}(x,y,z)=\frac{1}{2}\int\limits_{0}^{\infty}\frac{dt}{\sqrt{(t+x)(t+y)(t+z)}}.
\end{equation}
It is symmetric in $x$, $y$ and $z$, homogenous in all variables
of degree $-1/2$ and related to $F(\lambda,k)$ by
\begin{equation}\label{eq:F-Rf}
F(\lambda,k)={\lambda}R_{F}(1-\lambda^2,1-k^2\lambda^2,1).
\end{equation}

Asymptotic formulas for $F(\lambda,k)$ near the point $(1,1)$
appeared in \cite{Byrd,CarlsonBook,Gustafson,Nellis}, but the
first complete asymptotic series with error bounds at each order
of approximation was given by B.C.\,Carlson and J.L.\,Gustafson in
terms of the first symmetric standard elliptic integral $R_{F}$ in
\cite{CG1}. As is clear from (\ref{eq:F-Rf}) and homogeneity, the
case $\lambda,k\to{1}$ for $F(\lambda,k)$ is equivalent to the
case $z\to{\infty}$ with bounded $x$ and $y$ for $R_{F}(x,y,z)$.
The first two approximations from \cite{CG1} translated into our
notation read:
\begin{equation}\label{eq:CG1}
F(\lambda,k)=\lambda\ln\frac{4}{\sqrt{1-\lambda^2}+\sqrt{1-k^2\lambda^2}}+\theta_1F(\lambda,k),
\end{equation}
with relative error bound
\begin{equation}\label{eq:CG1-error}
\frac{(2-\lambda^2(1+k^2))\ln(1-k^2\lambda^2)}{4\ln[(1-k^2\lambda^2)/16]}<\theta_1<\frac{2-\lambda^2(1+k^2)}{4},
\end{equation}
and
\begin{multline}\label{eq:CG2}
F(\lambda,k)=\frac{\lambda}{4}\left[(6-\lambda^2(1+k^2))\ln\frac{4}{\sqrt{1-\lambda^2}+\sqrt{1-k^2\lambda^2}}
-2+\lambda^2(1+k^2)+\sqrt{(1-\lambda^2)(1-k^2\lambda^2)}\right]
\\+\theta_2F(\lambda,k),
\end{multline}
with relative error bound
\begin{equation}\label{eq:CG2-error}
\frac{9(1-k^2\lambda^2)^2\ln(1-k^2\lambda^2)}{64\ln[(1-k^2\lambda^2)/16]}<\theta_2<\frac{3(1-k^2\lambda^2)^2}{8}.
\end{equation}
The authors also provide more precise approximations at the price
of having the first complete elliptic integral in them:
\begin{equation}\label{eq:CG3}
F(\lambda,k)=\frac{2}{\pi}K\!\!\left(\sqrt{1-k^2}\right)\ln\frac{4}{\sqrt{1-\lambda^2}+\sqrt{1-k^2\lambda^2}}-\delta_1
\end{equation}
\begin{equation}\label{eq:CG4}
=\frac{2}{\pi}K\!\!\left(\sqrt{1-k^2}\right)\ln\frac{4}{\sqrt{1-\lambda^2}+\sqrt{1-k^2\lambda^2}}
-\frac{1}{4}\left(2-\lambda^2-k^2\lambda^2-\sqrt{(1-\lambda^2)(1-k^2)}\right)+\delta_2,
\end{equation}
where absolute errors have bounds given by
\begin{equation}\label{eq:CG34-error}
\frac{1-k^2\lambda^2}{8}<\delta_1<\frac{(1-k^2\lambda^2)\ln(4)}{k^2\lambda^2},~~~
\frac{9(1-k^2\lambda^2)^2}{64}<\delta_2<\frac{3(1-k^2\lambda^2)^2\ln(2)}{2k^2\lambda^2}.
\end{equation}
The problem of finding complete asymptotic expansion for
$z\to{\infty}$ with bounded $x$ and $y$ solved by Carlson and
Gustafson in \cite{CG1} for the first symmetric standard elliptic
integral $R_{F}$ has been solved for the other types of symmetric
standard EIs by J.L.\,L\'{o}pez in \cite{Lopez}. This research has
been continued in \cite{Lopez1}, where complete asymptotic
expansions are found for all symmetric standard elliptic integrals
when two variables tend to infinity thus settling the question in
principle for symmetric standard EIs. The methods used in the
above papers are either based on Mellin transform technique
\cite{CG1} or distributional approach \cite{Lopez,Lopez1}.  Recent
advances in Mellin transform technique can be found in
\cite{Lopez2}.

As one clearly sees from the error bounds (\ref{eq:CG1-error}),
(\ref{eq:CG2-error}) and (\ref{eq:CG34-error}) the corresponding
approximations are only asymptotic when \emph{both} variables
$\lambda$ and $k$ approach one.  The approximations derived in
this paper are of somewhat different character in that only one of
the variables needs to approach one while the other is allowed to
behave arbitrarily including approaching one as well. Hence,
Theorems~\ref{th:main1} and \ref{th:main2} below provide
asymptotic approximations for $F(\lambda,k)$ of any order for
$\lambda$ and $k$ lying on any curve having the endpoint at the
side $\lambda=1$ or  $k=1$ (including the logarithmic singularity
$\lambda=k=1$) of the unit square $[0,1]\times[0,1]$ in the
$(\lambda,k)$-plane. The coefficients of our first expansion
(\ref{eq:bottomexp}) are expressed recursively in terms of
elementary function.  The second expansion
(\ref{eq:Fasymptop-final}) contains the first complete elliptic
integral minus an elementary function also computed recursively.
Each approximation is accompanied by a two-sided error bound.  Our
derivation is based on simple rearrangements of certain
modifications of some known double series expansions,
hypergeometric summation algorithms and inequalities for
hypergeometric functions (some known and some new).  The resulting
approximations are very precise which is demonstrated in the last
section of the paper containing numerical examples and a
comparison with (\ref{eq:CG1}) and (\ref{eq:CG2}).

Expansion (\ref{eq:Fasymptop-final}) may be combined with
asymptotically precise inequalities for the first complete
elliptic integral found in \cite{Sitnik1}.  These inequalities can
be further improved by employing integral representations and
using the method of refining the Cauchy-Bunyakowsky integral
inequality developed in \cite{Sitnik2} - \cite{Sitnik4}.
Computations show good precision of these results near the
singularity.

\paragraph{2. Expansions of B.\,Radon and R.\,Kelisky revisited.}
In this section we derive two auxiliary expansions which will
serve as starting points for our main results formulated in
sections~3 and 4. The first expansion can be viewed as a different
guise of a known expansion due to Brigitte Radon \cite{Radon},
while the second follows from an expansion due to Richard Kelisky
\cite{Kelisky} by some hypergeometric transformations. The error
bounds found in this section appear to be new.

To keep the exposition as self-contained as possible we will give
an independent derivation of a modified Radon's expansion.  To
this end we need the following lemma.
\begin{lemma}\label{lm:Gauss-exp1}
For an integer $j\geq{0}$ the following identity is true\emph{:}
\begin{equation}\label{eq:Gauss-exp1}
\int\limits_{0}^{\lambda}\frac{t^{2j}dt}{(1-t^2)^{j+1}}=(-1)^j\frac{(1/2)_j}{2j!}\ln\left(\frac{1+\lambda}{1-\lambda}\right)+
\frac{1}{2j\lambda}\sum\limits_{n=0}^{j-1}(-1)^{n}
\frac{(1/2-j)_n}{(1-j)_n}\left(\frac{\lambda^2}{1-\lambda^2}\right)^{j-n},
\end{equation}
where the second term is zero for $j=0$. Here
\[
(a)_n=a(a+1)(a+2)\cdots(a+n-1),~~~~(a)_0=1,
\]
is the Pochhammer symbol \emph{(}or shifted factorial\emph{)}.
\end{lemma}

\textbf{Proof.} Euler's integral representation for the Gauss
hypergeometric function reads after a simple variable change:
\begin{equation}\label{eq:ourEuler}
\int\limits_{0}^{\lambda}\frac{t^{2j}dt}{(1-t^2)^{n}}=\frac{\lambda^{2j+1}}{2j+1}
{_{2}F_{1}}(n,j+1/2;j+3/2;\lambda^2).
\end{equation}
Writing (\ref{eq:ourEuler}) for $n=j+1$, employing differentiation
rule \cite[formula 2.1(7)]{Bat1} and representation \cite[formula
2.8(14)]{Bat1} for ${_{2}F_{1}}(1,1/2;3/2;\lambda^2)$ we compute:
\[
\int\limits_{0}^{\lambda}\frac{t^{2j}dt}{(1-t^2)^{j+1}}=\frac{\lambda^{2j+1}}{2j+1}{_{2}F_{1}}(j+1,j+1/2;j+3/2;\lambda^2)
=\frac{\lambda^{2j+1}}{j!}
\left(\frac{d}{du}\right)^j\left[\frac{\arctanh(\sqrt{u})}{\sqrt{u}}\right]_{u=\lambda^2}
\]\[
=(-1)^j\frac{(1/2)_j}{2j!}\ln\left(\frac{1+\lambda}{1-\lambda}\right)+
\frac{1}{2\lambda}\sum\limits_{k=0}^{j-1}\sum\limits_{m=0}^{j-1-k}\frac{(-1)^{k+m}(1/2)_k(1/2)_m}{(j-k)k!m!}
\left(\frac{\lambda^2}{1-\lambda^2}\right)^{j-k-m}
\]
\[
=(-1)^j\frac{(1/2)_j}{2j!}\ln\left(\frac{1+\lambda}{1-\lambda}\right)+\left(\frac{\lambda^2}{1-\lambda^2}\right)^j
\frac{1}{2\lambda}\sum\limits_{n=0}^{j-1}\left(\frac{1-\lambda^2}{\lambda^2}\right)^{n}(-1)^{n}
\sum\limits_{k=0}^{n}\frac{(1/2)_k(1/2)_{n-k}}{(j-k)k!(n-k)!}
\]\[
=(-1)^j\frac{(1/2)_j}{2j!}\ln\left(\frac{1+\lambda}{1-\lambda}\right)+
\frac{1}{2j\lambda}\sum\limits_{n=0}^{j-1}(-1)^{n}
\frac{(1/2-j)_n}{(1-j)_n}\left(\frac{\lambda^2}{1-\lambda^2}\right)^{j-n}.~~~~\square
\]

\begin{theo}\label{th:F-first}
For $\lambda$ and $k$  satisfying
\begin{equation}\label{eq:k-lambda-ineq}
1-k^2<(1-\lambda^2)/\lambda^2,
\end{equation}
and an integer $N\geq{0}$ the expansion
\begin{multline}\label{eq:F-Radon}
F(\lambda,k)=\frac{1}{2}\ln\frac{1+\lambda}{1-\lambda}\sum\limits_{j=0}^{N}\frac{(1/2)_j(1/2)_j}{(j!)^2}(1-k^2)^{j}
\\
+\frac{1}{2\lambda}\sum\limits_{j=1}^{N}\sum\limits_{n=0}^{j-1}
(-1)^{j+n}\frac{(1/2)_j(1/2-j)_n}{jj!(1-j)_n}(1-k^2)^{j}\left(\frac{\lambda^2}{1-\lambda^2}\right)^{j-n}+R_{1,N}(\lambda,k)
\end{multline}
holds true.  The bound for the remainder is given by
\begin{equation}\label{eq:R1Fest}
|R_{1,N}(\lambda,k)|\leq
\frac{\lambda(1/2)_{N+1}}{2(N+1)(N+1)!}\left[\frac{\lambda^2(1-k^2)}{1-\lambda^2}\right]^{N+1}.
\end{equation}
\end{theo}
\rem  It is clear from the error bound (\ref{eq:R1Fest}) that
expansion (\ref{eq:F-Radon}) is asymptotic for
$(1-k)/(1-\lambda)\to{0}$ and convergent for $\lambda$ and $k$
satisfying (\ref{eq:k-lambda-ineq}).

\textbf{Proof.} Put ${k^\prime}^2=1-k^2$. Expanding
$\left[1+({k^\prime}^2t^2)/(1-t^2)\right]^{-1/2}$ into the
binomial series and  interchanging summation and integration we
compute:
\[
F(\lambda,k)=\int\limits_{0}^{\lambda}\frac{dt}{\sqrt{(1-t^2)(1-t^2+{k^\prime}^2t^2)}}=
\int\limits_{0}^{\lambda}\frac{dt}{1-t^2}\left(1+\frac{{k^\prime}^2t^2}{1-t^2}\right)^{-1/2}
\]\[
=\int\limits_{0}^{\lambda}\frac{dt}{1-t^2}
\left(\sum\limits_{j=0}^{\infty}(-1)^j\frac{(1/2)_j}{j!}\frac{{k^\prime}^{2j}t^{2j}}{(1-t^2)^j}\right)
\]\[
=
\sum\limits_{j=0}^{N}(-1)^j\frac{(1/2)_j}{j!}{k^\prime}^{2j}\int\limits_{0}^{\lambda}\frac{t^{2j}dt}{(1-t^2)^{j+1}}
+
\sum\limits_{j=N+1}^{\infty}(-1)^j\frac{(1/2)_j}{j!}{k^\prime}^{2j}\int\limits_{0}^{\lambda}\frac{t^{2j}dt}{(1-t^2)^{j+1}}.
\]
Writing the first integral on the right-hand side as
(\ref{eq:Gauss-exp1}), we obtain (\ref{eq:F-Radon}) with $R_{1,N}$
given by
\[
R_{1,N}(\lambda,k)=\sum\limits_{j=N+1}^{\infty}(-1)^j\frac{(1/2)_j}{j!}(1-k^2)^{j}\int\limits_{0}^{\lambda}\frac{t^{2j}dt}{(1-t^2)^{j+1}}.
\]
This series is obviously alternating.  The following estimate
shows that each term is smaller in absolute value than the
previous one:
\[
\frac{(1/2)_{j+1}}{(j+1)!}(1-k^2)^{j+1}\int\limits_{0}^{\lambda}\frac{t^{2j+2}dt}{(1-t^2)^{j+2}}=
\frac{(1/2)_j(1/2+j)}{j!(j+1)}(1-k^2)^{j}\int\limits_{0}^{\lambda}\frac{t^{2j}}{(1-t^2)^{j+1}}\frac{t^2(1-k^2)}{(1-t^2)}dt
\]\[
\leq\frac{(1/2)_j}{j!}(1-k^2)^{j}\int\limits_{0}^{\lambda}\frac{t^{2j}}{(1-t^2)^{j+1}}\frac{\lambda^2(1-k^2)}{(1-\lambda^2)}dt\leq
\frac{(1/2)_j}{j!}(1-k^2)^{j}\int\limits_{0}^{\lambda}\frac{t^{2j}dt}{(1-t^2)^{j+1}}.
\]
The last inequality is due to (\ref{eq:k-lambda-ineq}).   Hence,
we are in the position to apply the Leibnitz convergence test
which implies that the remainder term $R_{1,N}(\lambda,k)$ does
not exceed
\[
\frac{(1/2)_{N+1}}{(N+1)!}(1-k^2)^{N+1}\int\limits_{0}^{\lambda}\frac{t^{2N+2}dt}{(1-t^2)^{N+2}}.
\]
We will prove the following asymptotically exact (as
$\lambda\to{1}$) estimate
\begin{equation}\label{eq:lambdaint_est}
f_1(\lambda)\equiv\int\limits_{0}^{\lambda}\frac{t^{2a}dt}{(1-t^2)^{a+1}}\leq\frac{\lambda^{2a+1}}{2a(1-\lambda^2)^a}\equiv{f_2(\lambda)}
\end{equation}
valid for all $\lambda\in(0,1)$ and $a>0$.  Indeed,
$f_1(0)=f_2(0)=0$ and
\[
\frac{f_1'(\lambda)}{f_2'(\lambda)}=\frac{2a}{2a+1-\lambda^2}<1,
~~~\lambda\in(0,1).
\]
The estimate (\ref{eq:lambdaint_est}) immediately leads to
(\ref{eq:R1Fest}).~$\square$

\rem One can verify that expansion (\ref{eq:F-Radon}) is a
different form of the expansion
\begin{equation}\label{eq:Radon}
F(\lambda,k)=\sum\limits_{j=0}^{\infty}\frac{(2j)!}{(j!)^3}\left[\frac{(1-k^2)\lambda}{8\sqrt{1-\lambda^2}}\right]^j
\mathrm{Q}^j_j\left(\frac{1}{\lambda}\right)
\end{equation}
due to B.\,Radon (see \cite{Radon}).  Here $\mathrm{Q}^j_j$
denotes the Legendre function of the second kind. Indeed, using
representation (\ref{eq:ourEuler}) and the formula
\[
{_{2}F_{1}}(j+1,j+1/2;j+3/2;\lambda^2)=\frac{(-1)^j(2j+1)}{2^jj!\lambda^{j+1}(1-\lambda^2)^{j/2}}\mathrm{Q}^j_j\left(\frac{1}{\lambda}\right)
\]
instead of (\ref{eq:Gauss-exp1}) in the proof of
Theorem~\ref{th:F-first} we can get (\ref{eq:Radon}).

The function ${_2F_1}(-n,1/2;1;x)$ creeps up into our
considerations on several occasions.  It can be expressed in terms
of Legendre polynomials via
\begin{equation}\label{eq:F-Legendre0}
{_2F_1}(-n,1/2;1;x)=(1-x)^{n/2}P_n\left(\frac{(1-x)^{1/2}+(1-x)^{-1/2}}{2}\right)=(1-x)^{n/2}P_n\left(\frac{2-x}{2\sqrt{1-x}}\right)
\end{equation}
(see \cite[formulas 7.3.1(175)]{Prud3}). Using the first Laplace
integral
\[
P_n(z)=\frac{1}{\pi}\int\limits_{0}^{\pi}(z+\sqrt{z^2-1}\cos\varphi)^nd\varphi
\]
for $P_n(z)$ (see \cite[formula (4.8.10)]{Szego})), we obtain:
\begin{equation}\label{eq:Lapl-2F1}
{_2F_1}(-n,1/2;1;x)=\frac{1}{\pi}\int\limits_{0}^{\pi}\left(1-x\sin^2\frac{\varphi}{2}\right)^nd\varphi.
\end{equation}
We summarize the required knowledge about this function in the
following lemma.

\begin{lemma}
\emph{a)} The function $F_n(x)={_2F_1}(-n,1/2;1;x)$ is monotone
decreasing for $x\in[0,1]$ and all non-negative integers $n$ with
bounds
\begin{equation}\label{eq:2F1bound1}
F_n(1)=\frac{(1/2)_n}{n!}\leq {_2F_1}(-n,1/2;1;x)\leq F_n(0)=1.
\end{equation}

\emph{b)} For $x\in[1,2]$ the function $F_n(x)$ is monotone
decreasing when $n$ is odd with bounds
\begin{equation}\label{eq:2F1bound2}
F_n(2)=0\leq {_2F_1}(-n,1/2;1;x)\leq
F_n(1)=\frac{(1/2)_n}{n!}\leq{1},
\end{equation}
and has a single minimum at $x_{\mathrm{min}}\in(1,2)$ when $n$ is
even with bounds
\begin{equation}\label{eq:2F1bound3}
0<{_2F_1}(-n,1/2;1;x)\leq F_n(2)=\frac{n!}{2^n(n/2)!^2}\leq{1}.
\end{equation}

\emph{c)} For $x>2$ the function $F_n(x)$ has the sign $(-1)^n$
and increasing (decreasing) for even (odd) $n$ with the bound
\begin{equation}\label{eq:2F1bound4}
|{_2F_1}(-n,1/2;1;x)|\leq(x-1)^n.
\end{equation}

\emph{d)} The following identity holds true
\begin{equation}\label{eq:2F1identity}
{_2F_1}(-n,1/2;1;1-x)=\frac{(1/2)_n}{n!}{_2F_1}(-n,1/2;1/2-n;x).
\end{equation}
\end{lemma}

\textbf{Proof.}  There are many ways to prove this lemma.  We
present a self-contained proof based on representation
(\ref{eq:Lapl-2F1}).

1. Let $0\leq{x}<1 $. Then (\ref{eq:Lapl-2F1}) shows that $F_n(x)$
is decreasing, so that $F_n(1)\leq F_n(x)\leq{F_n(0)}$. Clearly
$F_n(0)=1$, while $F_n(1)=(1/2)_n/n!$ is the celebrated
Chu-Vandermonde identity and a) follows.

2. Let $1\leq{x}\leq{2}$, $n=2k+1$ and $k\geq0$ is an integer.
Since $1-x\sin^2\frac{\varphi}{2}$ is decreasing in $x$,  we infer
from (\ref{eq:Lapl-2F1}) that $F_n(x)$ is also decreasing. Hence,
$F_n(2)\leq{F_n(x)}\leq F_n(1)$ and
\begin{equation*}
F_{2k+1}(2)=\frac{1}{\pi}\int_0^{\pi}\cos^{2k+1}\phi\,d\phi=0
\end{equation*}
which proves (\ref{eq:2F1bound2}).

3. Let $1\leq{x}\leq{2}$, $n=2k$ and $k\geq{0}$ is an integer. It
is then obvious from (\ref{eq:Lapl-2F1}) that $F''_n(x)>0$, so
that $F'_n(x)$ is increasing. It is also clear that $F'_n(1)<0$.
On the other endpoint we have
\begin{equation*}
F'_{2k}(2)=-\frac{2k}{\pi}\int_0^{\pi}\cos^{2k-1}\phi\frac{1-\cos\phi}{2}
\,d\phi=0+\frac{k}{\pi}\int_0^{\pi}\cos^{2k}\phi\,d\phi>0.
\end{equation*}
Monotonicity of $F'_{2k}(x)$ implies that there is a single
minimum at a point  $x_{\min}\!\in\!(1,2)$ and hence
$0<F_{2k}(x)<\max(F_{2k}(1),F_{2k}(2))$. The value of $F_{2k}(2)$
can be computed from (\ref{eq:Lapl-2F1}) or found in \cite[formula
7.3.8(2)]{Prud3}. To prove (\ref{eq:2F1bound3}) we need to show
that $F_{2k}(1)<F_{2k}(2)$ or
\begin{equation}
\frac{(\frac{1}{2})_n}{n!}<\frac{n!}{2^n(\frac{n}{2}!)^2}~\Leftrightarrow~
\frac{\Gamma(2k+1/2)}{\Gamma(2k+1)}<\frac{\Gamma(k+1/2)}{\Gamma(k+1)}.
\end{equation}
The last inequality is clearly true for $k=1$. Its validity for
any integer $k$ follows by induction on $k$, where the step of
induction is secured by the elementary inequality
\begin{equation*}
\frac{(2k+3/2)(2k+1/2)}{(2k+2)(2k+1)}<\frac{k+1/2}{k+1}.
\end{equation*}
Moreover,
\begin{equation*}
\frac{n!}{2^n(\frac{n}{2}!)^2}=\frac{\Gamma(k+1/2)}{\Gamma(k+1}<1~\Rightarrow~\text{b)}.
\end{equation*}

4. Suppose $x>2$, $n=2k$. From (\ref{eq:Lapl-2F1}) we see that
$F_{2k}(x)>0$, $F''_{2k}(x)>0$ and $F'_{2k}(2)>0$ which implies
that $F'_{2k}(x)>0$ and $F_{2k}(x)$ is increasing.

5. Suppose $x>2$, $n=2k+1$.  As before $F_{2k+1}(2)=0$,
$F'_{2k+1}(x)<0$ and consequently  $F_{2k+1}(x)<0$ and is
decreasing.

6. From the elementary inequality
\begin{equation*}
|1-\epsilon x|\leq|1-x|,~~x\geq2,~~0\leq\epsilon\leq{1}
\end{equation*}
and (\ref{eq:Lapl-2F1}) we get (\ref{eq:2F1bound4}) by choosing
$\epsilon=\sin^2\frac{\phi}{2}$.

7. Finally, identity (\ref{eq:2F1identity}) is the limiting case
of the well-known analytic extension formula for ${_2F_1}$, see
\cite[formula 15.3.6]{AS}. This identity can also be proved by
writing $F_n$ as the Legendre polynomial as in
(\ref{eq:F-Legendre0}) and applying  \cite[formulas
7.3.1(175)-(176)]{Prud3}.~$\square$

\begin{theo}\label{th:Kelisky-modified} For $\lambda$, $k$ satisfying
\begin{equation}\label{eq:region3}
(1-k^2)/k^2>1-\lambda^2
\end{equation}
and a positive integer $N$ the following expansion holds
true\emph{:}
\begin{equation}\label{eq:Freverse-exp2}
F(\lambda,k)=K(k)-\left[\frac{1-\lambda^2}{1-k^2}\right]^{1/2}
\sum\limits_{m=0}^{N-1}\frac{(1-\lambda^2)^{m}}{2m+1}
{_{2}F_{1}}(-m,1/2;1;(1-k^2)^{-1})+R_{2,N}(\lambda,k),
\end{equation}
where $K(k)$ is the complete elliptic integral of the first kind.
The bound for the remainder term is given by
\begin{equation}\label{eq:R2Fest}
|R_{2,N}(\lambda,k)|\leq \frac{1}{2N+1}
\left[\frac{(1-\lambda^2)k^2}{(1-k^2)}\right]^N\frac{\sqrt{(1-\lambda^2)(1-k^2)}}{1+k^2\lambda^2-2k^2}
\end{equation}
for $1/2\leq{k^2}<1$, and
\begin{equation}\label{eq:R2Fest1}
|R_{2,N}(\lambda,k)|\leq\frac{\lambda^{-2}}{2N+1}(1-\lambda^2)^N\sqrt{\frac{1-\lambda^2}{1-k^2}}
\end{equation}
for $0<{k^2}\leq{1/2}$.
\end{theo}
\rem  It is clear from the error bound
(\ref{eq:R2Fest})-(\ref{eq:R2Fest1}) that expansion
(\ref{eq:Freverse-exp2}) is asymptotic for
$(1-\lambda)/(1-k)\to{0}$ and convergent for $\lambda$ and $k$
satisfying (\ref{eq:region3}).

\textbf{Proof.} We begin with an expansion around $\lambda=k=0$
given by R.P.~Kelisky in \cite{Kelisky} (see also
\cite{Carlson1}):
\[
F(\lambda,k)=\sum\limits_{m=0}^{\infty}\frac{\lambda^{2m+1}}{(2m+1)}\frac{(1/2)_m}{m!}{_2F_1}(-m,1/2;1/2-m;k^2).
\]
An application of the identity (\ref{eq:2F1identity}) transforms
it into the expansion
\begin{equation}\label{eq:F-Gauss}
F(\lambda,k)=\sum\limits_{m=0}^{\infty}\frac{\lambda^{2m+1}}{2m+1}{_2F_1}(-m,1/2;1;1-k^2),
\end{equation}
valid for fixed $0<\lambda<1$ and $|k|<1/\lambda$. To make the
next  step we need the reflection-type relation
\begin{equation}\label{eq:F-relation}
F(\lambda,k)=K(k)-\frac{1}{\sqrt{1-k^2}}F\left(\sqrt{1-\lambda^2},\sqrt{-k^2/(1-k^2)}\right),
\end{equation}
which can be easily verified by representing the integral over
$(0,\lambda)$ in (\ref{eq:F-defined}) as the difference of
integrals over $(0,1)$ and $(\lambda,1)$ and introducing the new
integration variable $u^2=1-t^2$. The branch of
$\sqrt{1-\lambda^2}$ is chosen so that it is positive for positive
values of $1-\lambda^2$. The branch choice of the second square
root is immaterial since $F$ depends on the squared second
argument only. Expanding the second term on the right-hand side of
(\ref{eq:F-relation}) into the series (\ref{eq:F-Gauss}) and
splitting the resulting series we get (\ref{eq:Freverse-exp2})
with the remainder given by
\begin{equation}\label{eq:R2-defined}
R_{2,N}(\lambda,k)=\left[\frac{1-\lambda^2}{1-k^2}\right]^{1/2}\sum\limits_{m=N}^{\infty}\frac{(1-\lambda^2)^{m}}{2m+1}
{_{2}F_{1}}(-m,1/2;1;(1-k^2)^{-1}).
\end{equation}
To obtain a bound for $R_{2,N}$ we invoke the estimate
(\ref{eq:2F1bound4}).  Substituting this estimate into
(\ref{eq:R2-defined}) yields:
\[
|R_{2,N}(\lambda,k)|\leq\left[\frac{1-\lambda^2}{1-k^2}\right]^{1/2}
\sum\limits_{m=N}^{\infty}\frac{1}{2m+1}\left[\frac{(1-\lambda^2)k^2}{(1-k^2)}\right]^m
\]
for $1/2\leq{k^2}<1$ and
\[
|R_{2,N}(\lambda,k)|\leq\left[\frac{1-\lambda^2}{1-k^2}\right]^{1/2}
\sum\limits_{m=N}^{\infty}\frac{(1-\lambda^2)^m}{2m+1}
\]
for $0<{k^2}\leq{1/2}$. Applying the inequality
\[
\sum\limits_{m=N}^{\infty}\frac{x^m}{2m+1}=x^N\sum\limits_{s=0}^{\infty}\frac{x^s}{2s+2N+1}\leq
\frac{x^N}{2N+1}\sum\limits_{s=0}^{\infty}x^s=\frac{x^{N}}{(2N+1)(1-x)}
\]
valid for $0<x<1$, we arrive at (\ref{eq:R2Fest}) and
(\ref{eq:R2Fest1}).~~$\square$

\paragraph{3. First asymptotic expansion.} Denote
\begin{equation}\label{eq:sn}
s_n(x)=\sum\limits_{j=n+1}^{\infty}\frac{(1/2)_j(1/2-j)_n}{j!j(1-j)_n}(-x)^j.
\end{equation}
The change of summation variable $m=j-n-1$ gives:
\begin{equation}\label{eq:sn-4F3}
s_n(x)=2\frac{[(1/2)_{n+1}]^2}{[(n+1)!]^2}(-x)^{n+1}
{_4F_3}\left(\begin{array}{c}1,1,3/2+n,3/2+n\\3/2,2+n,2+n\end{array}\vline-x\right).
\end{equation}
We derive a recurrence formula for $s_n(x)$ in terms of elementary
functions in the following lemma.
\begin{lemma}\label{lm:sn}
The functions $s_n(x)$ satisfy the  four-term  recurrence relation
\begin{equation}\label{eq:recurrence}
4(n+3)^2s_{n+3}=a_n(x)s_{n+2}(x)+b_n(x)s_{n+1}(x)+c_n(x)s_{n}(x)+h_n(x),
\end{equation}
where
\[
a_n(x)=8n^2+36n+42-x(2n+5)^2,
\]\[
b_n(x)=2x(4n^2+14n+13)-(2n+3)^2,
\]\[
c_n(x)=-4x(n+1)^2,
\]\[
h_n(x)=\frac{x(2n+5)(2n+3)^2+(n+3)(8n^2+24n+17)}{8(n+3)[(n+2)!]^2}[(3/2)_n]^2(-x)^{n+2},
\]
and the starting values for the recursion are given by
\begin{equation}\label{eq:s0}
s_0(x)=-2\ln\frac{1+\sqrt{1+x}}{2},
\end{equation}
\begin{equation}\label{eq:s1}
s_1(x)=\left(\frac{x}{2}-1\right)\ln\frac{1+\sqrt{1+x}}{2}-\frac{1}{2}\sqrt{1+x}+\frac{1}{2}+\frac{x}{2},
\end{equation}
\begin{equation}\label{eq:s2}
s_2(x)=\left(-\frac{9}{32}x^2+\frac{x}{4}-\frac{3}{4}\right)\ln\frac{1+\sqrt{1+x}}{2}+\left(\frac{9}{32}x-\frac{7}{16}\right)\sqrt{1+x}
+\frac{7}{16}+\frac{1}{8}x-\frac{21}{64}x^2.
\end{equation}
\end{lemma}

\textbf{Proof.} Relation (\ref{eq:recurrence}) can be proved by a
careful application of Sister Celine's or Zeilberger's algorithm
\cite{Koepf}.  Denote the generic term in (\ref{eq:sn}) by
\[
g(n,j)=\frac{(1/2)_j(1/2-j)_n}{j!j(1-j)_n}(-x)^j=2\frac{[(1/2)_j]^2(j-n-1)!}{[j!]^2(3/2)_{j-n+1}}(-x)^j.
\]
The $j$-free recurrence relation
\begin{multline}\label{eq:j-free}
(-8n^2-36n-42)g(n+2,j+1)+4(n+3)^2g(n+3,j+1)+x(2n+5)^2g(n+2,j)\\
+4x(n+1)^2g(n,k)-2x(4n^2+14n+13)g(n+1,j)+(2n+3)^2g(n+1,j+1)=0
\end{multline}
can be verified by a direct substitution.  The difference from the
standard algorithms for hypergeometric summation comes from the
fact that we have non-standard bounds for the summation index
(which should be over all integers for standard algorithms).  We
can, however, remedy this by noting that
\[
\sum\limits_{j=n+3}^{\infty}g(n+2,j+1)=s_{n+2}(x)-g(n+2,n+3),
\]\[
\sum\limits_{j=n+3}^{\infty}g(n+3,j+1)=s_{n+3}(x),
\]\[
\sum\limits_{j=n+3}^{\infty}g(n+2,j)=s_{n+2}(x),
\]\[
\sum\limits_{j=n+3}^{\infty}g(n,j)=s_n(x)-g(n,n+1)-g(n,n+2),
\]\[
\sum\limits_{j=n+3}^{\infty}g(n+1,j)=s_{n+1}(x)-g(n+1,n+2),
\]\[
\sum\limits_{j=n+3}^{\infty}g(n+1,j+1)=s_{n+1}-g(n+1,n+2)-g(n+1,n+3),
\]
and summing up the $j$-free recurrence (\ref{eq:j-free}) over the
range $j=n+3,n+4,\ldots$. Together with definition of $g(n,j)$
this yields (\ref{eq:recurrence}).

To evaluate the initial term
\[
s_0(x)=\sum\limits_{j=1}^{\infty}\frac{(1/2)_j}{jj!}(-x)^j
\]
in a closed form note the identities
\[
\frac{1}{x}\left[\frac{1}{\sqrt{1-x}}-1\right]=\frac{1}{x}\left[\sum\limits_{j=0}^{\infty}\frac{(1/2)_j}{j!}x^j-1\right]=\sum\limits_{j=1}^{\infty}\frac{(1/2)_j}{j!}x^{j-1},
\]
\begin{equation}\label{eq:zerosum}
\int\limits_{0}^{x}\left[\frac{1}{\sqrt{1-t}}-1\right]\frac{dt}{t}=\sum\limits_{j=1}^{\infty}\frac{(1/2)_j}{jj!}x^j=2\ln\frac{2}{1+\sqrt{1-x}},
\end{equation}
and substitute $x$ with $-x$. The expression (\ref{eq:s1}) for the
next term
\begin{equation}\label{eq:onesum}
s_1(x)=\sum\limits_{j=2}^{\infty}\frac{(1/2)_j(1/2-j)}{(1-j)jj!}(-x)^{j}
\end{equation}
is derived from the identity
\[
\sum\limits_{j=2}^{\infty}\frac{(1/2)_j(1/2-j)}{(1-j)jj!}x^j=
\sum\limits_{j=2}^{\infty}\frac{(1/2)_j}{jj!}x^{j}+\frac{1}{2}\sum\limits_{j=2}^{\infty}\frac{(1/2)_j}{(j-1)jj!}x^{j},
\]
and evaluations
\[
\sum\limits_{j=2}^{\infty}\frac{(1/2)_j}{jj!}x^{j}=\sum\limits_{j=1}^{\infty}\frac{(1/2)_j}{jj!}x^{j}-\frac{x}{2}
=2\ln\frac{2}{1+\sqrt{1-x}}-\frac{x}{2},
\]\[
\sum\limits_{j=2}^{\infty}\frac{(1/2)_j}{jj!}x^{j-2}=\frac{2}{x^2}\ln\frac{2}{1+\sqrt{1-x}}-\frac{1}{2x}
\]
both deduced from (\ref{eq:zerosum}).  Hence,
\[
\sum\limits_{j=2}^{\infty}\frac{(1/2)_j}{(j-1)jj!}x^{j-1}
=\int\limits_{0}^{x}\left(\frac{2}{t^2}\ln\frac{2}{1+\sqrt{1-t}}-\frac{1}{2t}\right)dt
=\left(1-\frac{2}{x}\right)\ln\frac{2}{1+\sqrt{1-x}}+\frac{1}{1+\sqrt{1-x}}
\]
and (\ref{eq:onesum}) follows.   Similar but more cumbersome
computations lead to formula (\ref{eq:s2}) for
$s_2(x)$.~~$\square$

The main result of this section is now formulated as follows.
\begin{theo}\label{th:main1}
For all $(\lambda,k)\in[0,1]\times[0,1]$ and an integer
$N\geq{1}$, the first elliptic integral admits the representation
\begin{equation}\label{eq:bottomexp}
F(\lambda,k)=\frac{1}{2}\ln\frac{1+\lambda}{1-\lambda}\sum\limits_{j=0}^{N}\frac{(1/2)_j(1/2)_j}{(j!)^2}(1-k^2)^{j}
+\frac{1}{2\lambda}\sum\limits_{n=0}^{N-1}\left(\frac{1-\lambda^2}{-\lambda^2}\right)^{\!n}
s_n\!\left(\frac{(1-k^2)\lambda^2}{1-\lambda^2}\right)+R_N(\lambda,k),
\end{equation}
where $s_n(\cdot)$ is found from
\emph{(\ref{eq:recurrence})}-\emph{(\ref{eq:s2})}. The remainder
term  is negative and satisfies
\begin{equation}\label{eq:bottomerror}
\frac{[(1/2)_{N+1}]^2(1-k^2)^{N}}{2[(N+1)!]^2}f_{N+1}(\lambda,k)\leq-R_N(\lambda,k)\leq\frac{[(1/2)_{N+1}]^2(1-k^2)^{N}}{2[(N+1)!]^2}f_N(\lambda,k),
\end{equation}
where the positive function
\begin{equation}\label{eq:fN}
f_N(\lambda,k)=\frac{1}{1-\alpha(1-k^2)}\left\{\frac{1}{\alpha\lambda\sqrt{1+\frac{1-\lambda^2}{\alpha\lambda^2(1-k^2)}}}
\ln\frac{\sqrt{1+\frac{1-\lambda^2}{\alpha\lambda^2(1-k^2)}}+1}{\sqrt{1+\frac{1-\lambda^2}{\alpha\lambda^2(1-k^2)}}-1}
-(1-k^2)\ln\frac{1+\lambda}{1-\lambda}\right\}_{\vert
\alpha=\frac{(N+1/2)^2}{(N+1)^2}}
\end{equation}
is bounded on every subset $E$ of the unit square, where
\begin{equation}\label{eq:klamba-bounded}
\sup\limits_{k,\lambda\in{E}}\frac{1-k}{1-\lambda}<\infty
\end{equation}
and is monotone decreasing in $N$.
\end{theo}

\rem  Error bound (\ref{eq:bottomerror}) shows that expansion
(\ref{eq:bottomexp}) is asymptotic as $k\to{1}$ along any curve
$E$ lying entirely inside the unit square with
(\ref{eq:klamba-bounded}) satisfied, including those with endpoint
$(1,1)$.  The  expansion is convergent for any fixed
$0<\lambda<1$, $0<k<1$.

\rem  If condition (\ref{eq:klamba-bounded}) is violated but
$(1-k)^m/(1-\lambda)$ remains bounded, then $m$-th and higher
approximations are still asymptotic.  In this case, however, it is
much more effective to use approximation
(\ref{eq:Fasymptop-final}) from Theorem~\ref{th:main2}.

\textbf{Proof}. For the incomplete elliptic integral of the first
kind we have
\begin{multline}\label{eq:F-asymp}
F(\lambda,k)=\frac{1}{2}\ln\frac{1+\lambda}{1-\lambda}\sum\limits_{j=0}^{\infty}\frac{(1/2)_j(1/2)_j}{(j!)^2}(1-k^2)^{j}
\\
+\frac{1}{2\lambda}\sum\limits_{j=1}^{\infty}\sum\limits_{n=0}^{j-1}
(-1)^{j+n}\frac{(1/2)_j(1/2-j)_n}{jj!(1-j)_n}(1-k^2)^{j}\left(\frac{\lambda^2}{1-\lambda^2}\right)^{j-n}
\end{multline}
according to (\ref{eq:F-Radon}).   Rearranging the double sum
according to the rule
\begin{equation}\label{eq:rearrange}
\sum\limits_{j=N+1}^{\infty}\sum\limits_{n=N}^{j-1}a_{n,j}=\sum\limits_{n=N}^{\infty}\sum\limits_{j=n+1}^{\infty}a_{n,j}
\end{equation}
(taken for $N=0$), we get
\begin{multline}\label{eq:F-asymp-new1}
F(\lambda,k)=\frac{1}{2}\ln\frac{1+\lambda}{1-\lambda}\sum\limits_{j=0}^{\infty}\frac{(1/2)_j(1/2)_j}{(j!)^2}(1-k^2)^{j}
\\
+\frac{1}{2\lambda}\sum\limits_{n=0}^{\infty}\sum\limits_{j=n+1}^{\infty}
(-1)^{j+n}\frac{(1/2)_j(1/2-j)_n}{jj!(1-j)_n}(1-k^2)^{j}\left(\frac{\lambda^2}{1-\lambda^2}\right)^{j-n}.
\end{multline}
Summing the first series for $j$ from $0$ to $N$ and the second
for $n$ from $0$ to $N-1$ and leaving the rest as a remainder we
obtain (\ref{eq:bottomexp}) by Lemma~\ref{lm:sn} and definition
(\ref{eq:sn}) of the functions $s_n$. The remainder term is thus
given by
\begin{multline*}
R_N(\lambda,k)=\frac{1}{2}\ln\frac{1+\lambda}{1-\lambda}\sum\limits_{j=N+1}^{\infty}\frac{(1/2)_j(1/2)_j}{(j!)^2}(1-k^2)^{j}
\\
+\frac{1}{2\lambda}\sum\limits_{n=N}^{\infty}\sum\limits_{j=n+1}^{\infty}
(-1)^{j+n}\frac{(1/2)_j(1/2-j)_n}{jj!(1-j)_n}(1-k^2)^{j}\left(\frac{\lambda^2}{1-\lambda^2}\right)^{j-n}.
\end{multline*}
To estimate $R_N$ we change the order of summations in the second
term according to the rule (\ref{eq:rearrange}) applied from right
to left.  This yields
\[
R_N(\lambda,k)=\frac{1}{2}\ln\frac{1+\lambda}{1-\lambda}\sum\limits_{j=N+1}^{\infty}\frac{(1/2)_j(1/2)_j}{(j!)^2}(1-k^2)^{j}
\]\[
+\frac{1}{2\lambda}\sum\limits_{j=N+1}^{\infty}(-1)^j\frac{(1/2)_j}{jj!}\left[\frac{\lambda^2(1-k^2)}{1-\lambda^2}\right]^{j}
\sum\limits_{n=N}^{j-1}\frac{(1/2-j)_n}{(1-j)_n}\left(\frac{1-\lambda^2}{-\lambda^2}\right)^{n}.
\]
Introducing the new summation variable $k=j-n$ and applying
standard hypergeometric summation algorithms as realized by Maple
''sum" procedure, we get
\begin{multline*}
\sum\limits_{n=N}^{j-1}\frac{(1/2-j)_n}{(1-j)_n}\left(\frac{1-\lambda^2}{-\lambda^2}\right)^{n}
=\frac{(1/2)_j}{(j-1)!}\left(\frac{1-\lambda^2}{-\lambda^2}\right)^{j}
\sum\limits_{k=1}^{j-N}\frac{(k-1)!}{(1/2)_k}\left(\frac{-\lambda^2}{1-\lambda^2}\right)^{k}
\\
=\frac{(1/2)_j}{(j-1)!}
\left\{2\left(\frac{1-\lambda^2}{-\lambda^2}\right)^{j-1}{_2F_1\left(1,1;3/2;\frac{-\lambda^2}{1-\lambda^2}\right)}\right.
\\
\left.-\frac{(j-N)!}{(1/2)_{j-N+1}}\left(\frac{1-\lambda^2}{-\lambda^2}\right)^{N-1}{_2F_1}\left(1,1+j-N;3/2+j-N;\frac{-\lambda^2}{1-\lambda^2}\right)
\right\}
\end{multline*}
for the inner sum.   Now use
\[
{_2F_1\left(1,1;3/2;\frac{-\lambda^2}{1-\lambda^2}\right)}=\frac{1-\lambda^2}{2\lambda}\ln\frac{1+\lambda}{1-\lambda}
\]
for the first term in braces and Euler's integral representation
\[
{_2F_1}\left(1,1+j-N;3/2+j-N;\frac{-\lambda^2}{1-\lambda^2}\right)=\frac{(1/2)_{j-N+1}}{(j-N)!}\int\limits_{0}^{1}\frac{t^{j-N}(1-t)^{-1/2}dt}{(1+t\lambda^2/(1-\lambda^2))}
\]\[
=-\frac{(1/2)_{j-N+1}}{(j-N)!}
\left(\frac{1-\lambda^2}{-\lambda^2}\right)^{j-N+1}
\int\limits_{0}^{\frac{\lambda^2}{1-\lambda^2}}\frac{(-u)^{j-N}du}{(1+u)\sqrt{1-u\frac{1-\lambda^2}{\lambda^2}}}
\]
for the second.  Substituting these expressions into the above
formula for $R_N(\lambda,k)$ and interchanging summation and
integration we arrive at
\begin{equation}\label{eq:Rn}
R_N(\lambda,k)=\frac{(-1)^N}{2\lambda}\int\limits_{0}^{\frac{\lambda^2}{1-\lambda^2}}
\frac{\left(1-u\frac{1-\lambda^2}{\lambda^2}\right)^{-\frac{1}{2}}}{u^N(1+u)}\,du
\sum\limits_{j=N+1}^{\infty}\!\!\frac{[(1/2)_j]^2}{(j!)^2}[-(1-k^2)u]^j.
\end{equation}
An easy computation shows that
\[
\sum\limits_{j=N+1}^{\infty}\!\!\frac{[(1/2)_j]^2}{(j!)^2}(-x)^j=\frac{(1/2)_{N+1}]^2}{[(N+1)!]^2}(-x)^{N+1}
{_3F_2}(1,N+3/2,N+3/2;N+2,N+2;-x).
\]
Hence,
\[
R_N(\lambda,k)=-\frac{(1-k^2)^{N+1}[(1/2)_{N+1}]^2}{2\lambda[(N+1)!]^2}\int\limits_{0}^{\frac{\lambda^2}{1-\lambda^2}}
\frac{{_3F_2}(1,N+3/2,N+3/2;N+2,N+2;-(1-k^2)u)udu}{(1+u)\left(1-u\frac{1-\lambda^2}{\lambda^2}\right)^{1/2}}.
\]
Next, we apply the inequality
\begin{equation}\label{eq:3F2ineq}
\frac{1}{1+\frac{(N+3/2)^2}{(N+2)^2}x}\leq{_3F_2}(N+3/2,N+3/2,1;N+2,N+2;-x)\leq\frac{1}{1+\frac{(N+1/2)^2}{(N+1)^2}x},
\end{equation}
valid for $x>0$.  A proof of this inequality will be given
elsewhere.  Thus we have $R_N(\lambda,k)<0$ and
\[
\frac{(1-k^2)^{N}[(1/2)_{N+1}]^2}{2[(N+1)!]^2}g(\alpha_2,\lambda,k)
\leq-R_N(\lambda,k)\leq\frac{(1-k^2)^{N}[(1/2)_{N+1}]^2}{2[(N+1)!]^2}g(\alpha_1,\lambda,k),
\]
where
\begin{multline}\label{eq:g-computed}
g(\alpha,\lambda,k)=\frac{1-k^2}{\lambda}\int\limits_{0}^{\frac{\lambda^2}{1-\lambda^2}}
\frac{udu}{\left[1+\alpha(1-k^2)u\right](1+u)\left(1-u\frac{1-\lambda^2}{\lambda^2}\right)^{1/2}}
\\
=\frac{1}{1-\alpha(1-k^2)}\left\{\frac{1}{\alpha\lambda\sqrt{1+\frac{1-\lambda^2}{\alpha\lambda^2(1-k^2)}}}
\ln\frac{\sqrt{1+\frac{1-\lambda^2}{\alpha\lambda^2(1-k^2)}}+1}{\sqrt{1+\frac{1-\lambda^2}{\alpha\lambda^2(1-k^2)}}-1}
-(1-k^2)\ln\frac{1+\lambda}{1-\lambda}\right\}
\end{multline}
and $\alpha_1=[(N+1/2)/(N+1)]^2$, $\alpha_2=[(N+3/2)/(N+2)]^2$.
Defining
\[
f_N(\lambda,k)=g([(N+1/2)/(N+1)]^2,\lambda,k)
\]
we obtain the error bound (\ref{eq:bottomerror}).  The statement
(\ref{eq:klamba-bounded}) about the boundedness of
$f_N(\lambda,k)$ follows from an examination of the right-hand
side of (\ref{eq:g-computed}). The monotonicity of
$f_N(\lambda,k)$ in $N$ is implied by monotonicity of
$g(\alpha,\lambda,k)$ in $\alpha$ which is clear from the integral
representation (\ref{eq:g-computed}). ~$\square$

\rem If simplicity is preferred to precision one can apply  the
elementary inequality
\[
{0}\leq{_3F_2}(1,N+3/2,N+3/2;N+2,N+2;-x)\leq{1}, ~~x\geq{0},
\]
instead of (\ref{eq:3F2ineq}).  Using this inequality and explicit
representation
\[
\int\limits_{0}^{\frac{\lambda^2}{1-\lambda^2}}
\frac{\left(1-u\frac{1-\lambda^2}{\lambda^2}\right)^{-\frac{1}{2}}udu}{(1+u)}
=\frac{2\lambda^2}{1-\lambda^2}-\lambda\ln\frac{1+\lambda}{1-\lambda},
\]
one gets the following error bound:
\[
{0}\leq-R_N(\lambda,k)\leq\frac{[(1/2)_{N+1}]^2(1-k^2)^{N+1}}{[(N+1)!]^2}
\left\{\frac{\lambda}{1-\lambda^2}-\frac{1}{2}\ln\frac{1+\lambda}{1-\lambda}\right\}.
\]
This error bound is relatively precise for $k>\lambda$ but loses
precision substantially for $\lambda>k$, while the bound
(\ref{eq:bottomerror}) is very precise for all values of
parameters.

The  first  order approximation obtained from (\ref{eq:bottomexp})
(see (\ref{eq:F1-bottom}) below),
\[
F_1(\lambda,k)=\frac{1}{2}\ln\frac{1+\lambda}{1-\lambda}+\frac{1}{\lambda}\ln\frac{2}{1+\sqrt{(1-\lambda^2k^2)/(1-\lambda^2)}}
+\frac{1-k^2}{8}\ln\frac{1+\lambda}{1-\lambda},
\]
has an amazing property to be correct asymptotic approximation for
$F(\lambda,k)$ not only as $k\to{1}$ but also as $\lambda\to{0}$
including the case when both $\lambda,k\to{0}$ along any curve.
Indeed one can easily check that
\[
F_1(\lambda,k)=\lambda+\left(\frac{25}{96}-\frac{1}{48}k^2+\frac{3}{32}k^4\right)\lambda^3+\O(\lambda^5)
\]
as $\lambda\to{0}$, while
\[
F(\lambda,k)=\lambda+\frac{1}{6}(1+k^2)\lambda^3+\O(\lambda^5),~~\lambda\to{0},
\]
so that
\[
F(\lambda,k)-F_1(\lambda,k)=\O(\lambda^3),~~\lambda\to{0}.
\]
Thus $F_1(\lambda,k)$ is a true approximation for two sides of the
unit square (including endpoints) - the side $\lambda=0$,
$k\in[0,1]$ and the side $k=1$, $\lambda\in[0,1]$. The same is
true for higher order approximations but the approximation order
for $\lambda\to{0}$ does not increase with $N$.

\paragraph{4. Second asymptotic expansion.}
Denote
\begin{equation}\label{eq:An}
A_n(x)=\sum\limits_{j=0}^{\infty}
\binom{n+j}{j}\frac{(-1)^j(1/2)_j}{(2(n+j)+1)j!}x^j.
\end{equation}
The following evident formula is more notational than meaningful:
\begin{equation}\label{eq:sum-found0}
A_n(x)=\frac{1}{2n+1}{_3F_2}\left(\frac{1}{2},n+\frac{1}{2},n+1;1,n+\frac{3}{2};-x\right).
\end{equation}
We give three representations for $A_n(x)$ in the following lemma.
The first is more convenient for computing explicit expressions,
the second is designed for easier estimation and the third
provides a hint for an alternative derivation of
Theorem~\ref{th:main2} (see details in Remark~\ref{rm:altmain2}).
\begin{lemma}\label{lm:An}
The following identities hold for the functions $A_n(x)$:
\begin{multline}\label{eq:An-represent}
A_n(x)=\frac{1}{n!}D^n_x\left[(-1)^n
\frac{(1/2)_n}{n!\sqrt{x}}\ln(\sqrt{1+x}+\sqrt{x})+\frac{\sqrt{1+x}}{2nx}\sum\limits_{j=0}^{n-1}(-1)^j\frac{(1/2-n)_j}{(1-n)_j}x^{n-j}\right]
\\
=\frac{1}{2x^{n+1/2}}\int\limits_{0}^{x}\frac{t^{n-1/2}}{\sqrt{1+t}}{_2F_1}\left(-n,1/2;1;\frac{t}{1+t}\right)dt
=\frac{1}{2x^{n+1/2}}\int\limits_{0}^{x}\frac{t^{n-1/2}}{(1+t)^{(n+1)/2}}P_n\left(\frac{2+t}{2\sqrt{1+t}}\right)dt,
\end{multline}
where the second term in brackets equals zero for $n=0$ and $D_x$
means differentiation in $x$.
\end{lemma}

\textbf{Proof}. From (\ref{eq:An}) $A_n(x)$ can be written  in the
form
\begin{equation}\label{eq:series-decomposed}
A_n(x)
=\frac{1}{n!}\sum\limits_{j=0}^{\infty}(j+1)(j+2)\ldots(j+n)\frac{1}{2n+2j+1}\frac{(-1)^j(1/2)_j}{j!}x^j.
\end{equation}
For a formal power series
\[
f(x)=\sum\limits_{j=0}^{\infty}a_jx^j
\]
we have ($D_x$ means differentiation in $x$)
\begin{equation}\label{eq:Dn}
D^n_xx^nf(x)=\sum\limits_{j=0}^{\infty}(j+1)(j+2)\ldots(j+n)a_jx^j
\end{equation}
and
\begin{equation}\label{eq:intf}
\frac{1}{2}x^{-n-1/2}\int\limits_{0}^{x}t^{n-1/2}f(t)dt
=\sum\limits_{j=0}^{\infty}\frac{a_j}{2n+2j+1}x^j.
\end{equation}
Putting $a_j=(-1)^j(1/2)_j/j!$ gives $f(x)=1/\sqrt{1+x}$ by the
binomial theorem. Combining (\ref{eq:Dn}) and (\ref{eq:intf}) we
obtain from (\ref{eq:series-decomposed}):
\begin{equation}\label{eq:sum-found1}
\sum\limits_{j=0}^{\infty}\binom{n+j}{j}\frac{(-1)^j(1/2)_j}{(2(n+j)+1)j!}x^j
=\frac{1}{2n!}D^n_xx^{-1/2}\int\limits_{0}^{x}\frac{t^{n-1/2}}{\sqrt{1+t}}dt.
\end{equation}
The integral on the right-hand side can be reduced to
(\ref{eq:Gauss-exp1}) by the variable change $y^2=t/(1+t)$,
$t=y^2/(1-y^2)$, $dt=2ydy/(1-y^2)^2$:
\begin{equation}\label{eq:varchange}
\int\limits_{0}^{x}\frac{t^{n-1/2}}{\sqrt{1+t}}dt=\int\limits_{0}^{x}t^{n-1}\sqrt{\frac{t}{1+t}}dt
=2\int\limits_{0}^{\sqrt{\frac{x}{1+x}}}\frac{y^{2n}dy}{(1-y^2)^{n+1}}
\end{equation}
\[
=2(-1)^n\frac{(1/2)_n}{n!}\ln(\sqrt{1+x}+\sqrt{x})+\frac{\sqrt{1+x}}{n\sqrt{x}}\sum\limits_{j=0}^{n-1}(-1)^j\frac{(1/2-n)_j}{(1-n)_j}x^{n-j}.
\]
Hence, from (\ref{eq:sum-found1}) and the above evaluation we
arrive at the first formula (\ref{eq:An-represent}).  An
alternative method of evaluating the right hand-side of
(\ref{eq:sum-found1}) is the following.  Make the variable change
$t=ux$:
\[
\int\limits_{0}^{x}\frac{t^{n-1/2}}{\sqrt{1+t}}dt=x^{n+1/2}\int\limits_{0}^{1}\frac{u^{n-1/2}du}{\sqrt{1+ux}},
\]
differentiate under the integral sign and apply the Leibnitz
formula
\[
D_x^nx^nf(x)=n!\sum\limits_{k=0}^{n}\binom{n}{k}\frac{x^k}{k!}D_x^kf(x)
\]
and the elementary formula
\[
D^k_x(1+ux)^{-1/2}=(-1)^k(1/2)_ku^k(1+ux)^{-k-1/2}
\]
to get
\[
\frac{1}{2n!}D^n_xx^{-1/2}\int\limits_{0}^{x}\frac{t^{n-1/2}}{\sqrt{1+t}}dt=
\frac{1}{2n!}\int\limits_{0}^{1}u^{n-1/2}\left[D^n_x\frac{x^{n}}{\sqrt{1+ux}}\right]du
\]\[
=\frac{1}{2n!}\int\limits_{0}^{1}u^{n-1/2}\left[n!\sum\limits_{k=0}^{n}\binom{n}{k}\frac{x^k}{k!}D^k_x\frac{1}{\sqrt{1+ux}}\right]du
\]\[
=\frac{1}{2n!}\int\limits_{0}^{1}u^{n-1/2}
\left[n!\sum\limits_{k=0}^{n}\binom{n}{k}\frac{x^k}{k!}(-1)^k(1/2)_ku^k(1+ux)^{-k-1/2}\right]du=
\]\[
=\frac{1}{2}\int\limits_{0}^{1}\frac{u^{n-1/2}}{\sqrt{1+ux}}{_2F_1}\left(-n,1/2;1;\frac{ux}{1+ux}\right)du.
\]
Finally, substituting back $t=ux$ we obtain:
\begin{equation}\label{eq:sum-found3}
A_n(x)=\frac{1}{2x^{n+1/2}}\int\limits_{0}^{x}\frac{t^{n-1/2}}{\sqrt{1+t}}{_2F_1}\left(-n,1/2;1;\frac{t}{1+t}\right)dt.
\end{equation}
The last equality in (\ref{eq:An-represent}) is a direct
consequence of (\ref{eq:F-Legendre0}).~~$\square$

\begin{theo}\label{th:main2}
For $(\lambda,k)\in[0,1]\times(0,1]$ and integer $N\geq{1}$, the
first elliptic integral admits the representation
\begin{equation}\label{eq:Fasymptop-final}
F(\lambda,k)=K(k)-
\frac{(1-\lambda^2)^{1/2}}{(1-k^2)^{1/2}}\sum\limits_{n=0}^{N-1}(1-\lambda^2)^{n}A_n\!\!\left(\frac{1-\lambda^2}{1-k^2}\right)+\tilde{R}_{N}(\lambda,k),
\end{equation}
where the functions $A_n(x)$ are found  from
\emph{(\ref{eq:An-represent})}. The remainder term is negative and
satisfies
\begin{multline}\label{eq:error-bound1}
\frac{(1-\lambda^2)^{N+1/2}}{2\lambda^2N\sqrt{2-\lambda^2-k^2}}\geq
-\tilde{R}_{N}(\lambda,k)\geq
\\
\frac{(1-\lambda^2)^{N-1/2}(1/2)_{N}}{2NN!}\!\left(\!\!
\sqrt{2-\lambda^2-k^2} -\frac{(1-k^2)}{2\sqrt{1-\lambda^2}}
\ln\!\left\{1+2\frac{1-\lambda^2}{1-k^2}+2\frac{\sqrt{(1-\lambda^2)(2-\lambda^2-k^2)}}{1-k^2}\right\}\right).
\end{multline}
\end{theo}
\rem The error bound (\ref{eq:error-bound1}) shows that the
expansion (\ref{eq:Fasymptop-final}) is asymptotic for
$\lambda\to{1}$ along any curve lying entirely inside the unit
square, including those with endpoint $(1,1)$. The expansion is
convergent for any fixed $0<\lambda<1$, $0<k<1$.

\textbf{Proof.} By Theorem~\ref{th:Kelisky-modified} for values of
$\lambda$ and $k$ satisfying (\ref{eq:region3}) we have expansion
(\ref{eq:Freverse-exp2}) which, after little modification,  can be
written as:
\begin{equation}\label{eq:Freverse-asymp}
F(\lambda,k)=K(k)-
\sum\limits_{m=0}^{\infty}\frac{(1-\lambda^2)^{m+1/2}}{2m+1}\sum\limits_{j=0}^{m}
\binom{m}{j}\frac{(-1)^j(1/2)_j}{j!(1-k^2)^{j+1/2}}.
\end{equation}
Changing the order of summation according to the rule
\[
\sum\limits_{m=0}^{\infty}\sum\limits_{j=0}^{m}a_{m,j}=\sum\limits_{n=0}^{\infty}\sum\limits_{i=n}^{\infty}a_{i,i-n}=
\sum\limits_{n=0}^{\infty}\sum\limits_{j=0}^{\infty}a_{n+j,j},
\]
we obtain the formula
\begin{equation}\label{eq:Fasymp-top}
F(\lambda,k)=K(k)-
\sum\limits_{n=0}^{\infty}(1-\lambda^2)^{n}\sum\limits_{j=0}^{\infty}
\frac{(1-\lambda^2)^{j+1/2}}{(1-k^2)^{j+1/2}}
\binom{n+j}{j}\frac{(-1)^j(1/2)_j}{(2(n+j)+1)j!},
\end{equation}
which in view of (\ref{eq:An}), can be split as follows:
\[
F(\lambda,k)=K(k)-
\frac{(1-\lambda^2)^{1/2}}{(1-k^2)^{1/2}}\sum\limits_{n=0}^{N-1}(1-\lambda^2)^{n}A_n\left(\frac{1-\lambda^2}{1-k^2}\right)+\tilde{R}_{N}(\lambda,k),
\]
with $A_n$ defined by (\ref{eq:An}) and
\[
\tilde{R}_{N}(\lambda,k)=-\frac{(1-\lambda^2)^{1/2}}{(1-k^2)^{1/2}}\sum\limits_{n=N}^{\infty}(1-\lambda^2)^{n}A_n\left(\frac{1-\lambda^2}{1-k^2}\right).
\]
As written in (\ref{eq:Fasymp-top}), the inner sum does not
converge unless $k<\lambda$. However, it was shown in
Lemma~\ref{lm:An} that $A_n(x)$ is an elementary function defined
for all $(\lambda,k)\in[0,1)\times[0,1)$. We will prove that the
outer sum converges for all such $\lambda$ and $k$. To this end
substitute the second formula (\ref{eq:An-represent}) for $A_n$
and apply the estimate from above from (\ref{eq:2F1bound1}) to
get:
\[
-\tilde{R}_{N}(\lambda,k)=\frac{1}{2}\sum\limits_{n=N}^{\infty}(1-k^2)^{n}
\!\!\!\!\!\!\!\!\!\int\limits_{0}^{(1-\lambda^2)/(1-k^2)}\!\!\!\!\!\!\!\!\!\frac{t^{n}}{\sqrt{t(1+t)}}{_2F_1}\left(-n,1/2;1;\frac{t}{1+t}\right)dt
\]\[
\leq\frac{1}{2}\!\!\!\!\!\!\!\!\!\int\limits_{0}^{(1-\lambda^2)/(1-k^2)}\!\!\!\!\!\frac{dt}{\sqrt{t(1+t)}}
\sum\limits_{n=N}^{\infty}[(1-k^2)t]^{n}
=\frac{(1-k^2)^{N}}{2}\!\!\int\limits_{0}^{(1-\lambda^2)/(1-k^2)}\!\!\!\!\!\frac{t^{N}dt}{\sqrt{t(1+t)}[1-(1-k^2)t]}
\]\[
\leq\frac{(1-k^2)^{N}}{2\lambda^2}\!\!\int\limits_{0}^{(1-\lambda^2)/(1-k^2)}\!\!\!\!\!\frac{t^{N}dt}{\sqrt{t(1+t)}}
=\frac{(1-k^2)^{N}}{\lambda^2}\int\limits_{0}^{\sqrt{\frac{1-\lambda^2}{2-k^2-\lambda^2}}}\frac{y^{2N}dy}{(1-y^2)^{N+1}},
\]
where we used (\ref{eq:varchange}) to obtain the last equality.
Now an application of inequality (\ref{eq:lambdaint_est}) with
$a=N$ gives the upper bound in (\ref{eq:error-bound1}).

To find a lower bound we again apply (\ref{eq:2F1bound1}) but this
time the estimate from below.  This yields:
\[
-\tilde{R}_{N}(\lambda,k)\geq\frac{1}{2}\int\limits_{0}^{(1-\lambda^2)/(1-k^2)}\frac{dt}{\sqrt{t(1+t)}}
\sum\limits_{n=N}^{\infty}\frac{(1/2)_n}{n!}(1-k^2)^nt^n
\]\[
=\frac{(1-k^2)^{N}(1/2)_{N}}{2N!}\!\!\!\!\!\int\limits_{0}^{(1-\lambda^2)/(1-k^2)}\!\!\!\!\!\frac{t^{N}}{\sqrt{t(1+t)}}
{_2F_1}(N+1/2,1;N+1;(1-k^2)t)dt.
\]
It follows from \cite[Theorem~1.10]{Ponnusamy} that
\[
{_2F_1}(N+1/2,1;N+1;y)\geq\frac{1}{\sqrt{1-y}}~~\text{for}~~~y\in(0,1).
\]
Consequently,
\[
-\tilde{R}_{N}(\lambda,k)\geq\frac{(1-k^2)^{N}(1/2)_{N}}{2N!}\!\!\!\!\!
\int\limits_{0}^{(1-\lambda^2)/(1-k^2)}\!\!\frac{t^{N-1/2}}{\sqrt{(1+t)(1-(1-k^2)t)}}dt
\]\[
\geq \frac{(1-k^2)^{N}(1/2)_{N}}{2N!}\!\!\!\!\!
\int\limits_{0}^{(1-\lambda^2)/(1-k^2)}\!\!t^{N-1}\sqrt{\frac{t}{1+t}}dt.
\]
Both functions $t^{N-1}$ and $\sqrt{\frac{t}{1+t}}$ are
non-decreasing for $t>0$, $N\geq{1}$.  Hence, we are in the
position to apply the Chebyshev inequality \cite[formula
IX(1.1)]{Mitrinovic} which results in:
\[
-\tilde{R}_{N}(\lambda,k)\geq
\frac{(1-k^2)^{N+1}(1/2)_{N}}{2(1-\lambda^2)N!} \!\!\!\!\!
\int\limits_{0}^{(1-\lambda^2)/(1-k^2)}\!\!t^{N-1}dt
\int\limits_{0}^{(1-\lambda^2)/(1-k^2)}\!\!\sqrt{\frac{t}{1+t}}dt
\]\[
=\frac{(1-\lambda^2)^{N-1/2}(1/2)_{N}}{2NN!}\!\left(\!\!
\sqrt{2-\lambda^2-k^2} -\frac{(1-k^2)}{2\sqrt{1-\lambda^2}}
\ln\!\left\{1+2\frac{1-\lambda^2}{1-k^2}+2\frac{\sqrt{(1-\lambda^2)(2-\lambda^2-k^2)}}{1-k^2}\right\}\right).
\]
$\square$

\rem\label{rm:altmain2} Expansion (\ref{eq:Fasymp-top})  can be
obtained directly from the definition of $F(\lambda,k)$ without a
use of Theorem~\ref{th:Kelisky-modified}. Indeed, substituting the
last formula (\ref{eq:An-represent}) into (\ref{eq:Fasymp-top}) we
get
\begin{multline}\label{eq:eq:Fasymp-top-int}
F(\lambda,k)=K(k)- \frac{1}{2}\sum\limits_{n=0}^{\infty}
(1-k^2)^{n}\int\limits_{0}^{\frac{1-\lambda^2}{1-k^2}}\frac{t^{n-1/2}}{(1+t)^{(n+1)/2}}P_n\left(\frac{2+t}{2\sqrt{1+t}}\right)dt
\\
=K(k)-
\frac{1}{2}\int\limits_{0}^{\frac{1-\lambda^2}{1-k^2}}\frac{dt}{\sqrt{t(1+t)}}
\sum\limits_{n=0}^{\infty}(1-k^2)^{n}\frac{t^{n}}{(1+t)^{n/2}}P_n\left(\frac{2+t}{2\sqrt{1+t}}\right).
\end{multline}
The generating function for the Legendre polynomials is given by
\cite{Szego}:
\[
\sum\limits_{n=0}^{\infty}z^nP_n(x)=\frac{1}{\sqrt{1-2xz+z^2}}.
\]
Hence,
\begin{equation}\label{eq:generating}
\sum\limits_{n=0}^{\infty}\frac{(1-k^2)^{n}t^{n}}{(1+t)^{n/2}}P_n\left(\frac{2+t}{2\sqrt{1+t}}\right)
=\frac{\sqrt{1+t}}{\sqrt{(1+k^2t)(1-t+k^2t)}},
\end{equation}
and so
\[
F(\lambda,k)=K(k)-\frac{1}{2}\int\limits_{0}^{\frac{1-\lambda^2}{1-k^2}}\frac{dt}{\sqrt{t(1+k^2t)(1-t+k^2t)}}.
\]
This formula can be obtained from the elementary relation
\[
F(\lambda,k)=K(k)-\int\limits_{\lambda}^{1}\frac{du}{\sqrt{(1-u^2)(1-k^2u^2)}}
\]
by the variable change $t=(1-u^2)/(1-k^2)$.  Thus the whole
process could be started from the above representation.  The
change of the integration variable and an application of
(\ref{eq:generating}) then give the expansion
(\ref{eq:eq:Fasymp-top-int}).  The representation of the general
term of this expansion as $n$-th derivative as in the first
formula (\ref{eq:An-represent}) can be then obtained by taking
$(2+t)/(2\sqrt{1+t})$ as a new integration variable in
(\ref{eq:eq:Fasymp-top-int}) and applying the Rodrigues formula
for the Legendre polynomials.

\paragraph{5. Results of computations.} In this section we present
several examples of computations with the expansions obtained
above. We also give a comparison with the approximations
(\ref{eq:CG1}) and (\ref{eq:CG2}) due to Carlson and Gustafson.

Consider expansion (\ref{eq:bottomexp}) first. From (\ref{eq:s0})
and (\ref{eq:s1}) we have:
\[
s_0\left(\frac{\lambda^2(1-k^2)}{1-\lambda^2}\right)
=2\ln\frac{2}{1+\sqrt{(1-\lambda^2k^2)/(1-\lambda^2)}},
\]
\[
s_1\left(\frac{\lambda^2(1-k^2)}{1-\lambda^2}\right)
=\left[1-\frac{\lambda^2(1-k^2)}{2(1-\lambda^2)}\right]\ln\frac{2}{1+\sqrt{(1-k^2\lambda^2)/(1-\lambda^2)}}
+\frac{1-k^2\lambda^2}{2(1-\lambda^2)}-\frac{1}{2}\sqrt{\frac{1-k^2\lambda^2}{1-\lambda^2}}.
\]
Hence, the first and the second order approximations read:
\begin{equation}\label{eq:F1-bottom}
F_1(\lambda,k)=\frac{1}{2}\ln\frac{1+\lambda}{1-\lambda}+\frac{1}{\lambda}\ln\frac{2}{1+\sqrt{(1-\lambda^2k^2)/(1-\lambda^2)}}
+\frac{1-k^2}{8}\ln\frac{1+\lambda}{1-\lambda},
\end{equation}
\begin{multline}\label{eq:F2-bottom}
F_2(\lambda,k)=F_1(k,\lambda)+\left(\frac{1-k^2}{4\lambda}-\frac{(1-\lambda^2)}{2\lambda^3}\right)\ln\frac{2}{1+\sqrt{(1-\lambda^2k^2)/(1-\lambda^2)}}
\\
-\frac{(1-k^2)\sqrt{1-k^2\lambda^2}}{4\lambda\sqrt{1-\lambda^2}+4\lambda\sqrt{1-k^2\lambda^2}}
+\frac{9}{128}(1-k^2)^2\ln\frac{1+\lambda}{1-\lambda}.
\end{multline}
Denote by $\Delta_N$ the difference between the upper and the
lower bounds in (\ref{eq:bottomerror}):
\begin{equation}\label{eq:DeltaN}
\Delta_N=\frac{[(1/2)_{N+1}]^2(1-k^2)^{N}}{2[(N+1)!]^2}(f_N(\lambda,k)-f_{N+1}(\lambda,k)).
\end{equation}
Approximation (\ref{eq:F1-bottom}) combined with inequality
(\ref{eq:bottomerror}) puts $F(\lambda,k)$ within an interval of
length $\Delta_1$, while (\ref{eq:F2-bottom}) puts $F(\lambda,k)$
within an interval of length $\Delta_2$.  Numerical results are
presented in Table~1.  The exact values of $F(\lambda,k)$ shown in
the tables below have been computed by Maple wth the required
number of precise digits guaranteed.
\renewcommand\arraycolsep{3pt}
\setlength{\arrayrulewidth}{0.2pt}%
\setlength{\doublerulesep}{0pt} {\small
$$
\begin{array}{!{\vrule width 1.2pt\relax}c|c|c!{\vrule width
1.2pt\relax}c|c|c!{\vrule width 1.2pt\relax}c|c|c!{\vrule width
1.2pt\relax}}
\hline\hline\hline\hline\hline &&&&&&&&\\[-10pt]
\lambda & k & F(\lambda,k) & \begin{array}{c}\text{1st order}\\[-3pt]\text{approx. (\ref{eq:F1-bottom})}\end{array} &
\begin{array}{c}\text{Absolute}\\[-3pt]\text{error}\end{array} &
\begin{array}{c}\text{Length of}\\[-3pt]\text{error range}\\[-3pt]\Delta_1\end{array} &
\begin{array}{c}\text{2nd order}\\[-3pt]\text{approx. (\ref{eq:F2-bottom})}\end{array}&
\begin{array}{c}\text{Absolute}\\[-3pt]\text{error}\end{array} &
\begin{array}{c}\text{Length of}\\[-3pt]\text{error range}\\[-3pt]\Delta_2\end{array}
\\[5pt]
\hline\hline\hline\hline\hline
.8 & .8   & 1.0178 & 1.0334 & -.01554 & .742\!\!\times\!\!10^{-3} & 1.0216 & -.00378 & .926\!\!\times\!\!10^{-4}\\
.9 & .9   & 1.3532 & 1.3652 & -.01198 & .657\!\!\times\!\!10^{-3} & 1.3547 & -.00153 & .427\!\!\times\!\!10^{-4}\\
.95 & .95 & 1.6861 & 1.6936 & -.00750 & .430\!\!\times\!\!10^{-3} & 1.6866 & -.4914\!\!\times\!\!10^{-3}& .143\!\!\times\!\!10^{-4}\\
.99 & .99 & 2.4708 & 2.4726 & -.00185 & .107\!\!\times\!\!10^{-3} & 2.4708 & -.2468\!\!\times\!\!10^{-4}& .721\!\!\times\!\!10^{-6}\\
.95 & .99 & 1.7951 & 1.7955 & -.405\!\!\times\!\!10^{-3} & .639\!\!\times\!\!10^{-5} & 1.7951 & -.554\!\!\times\!\!10^{-5}& .463\!\!\times\!\!10^{-7}\\
.99 & .999& 2.6240 & 2.6240 & -.253\!\!\times\!\!10^{-4} & .213\!\!\times\!\!10^{-6} & 2.6240 & -.350\!\!\times\!\!10^{-7}& .157\!\!\times\!\!10^{-9}\\
\hline\hline\hline\hline\hline
\end{array}
$$
}\textbf{Table~1.} \textsl{Numerical examples for approximations
(\ref{eq:F1-bottom}) and (\ref{eq:F2-bottom}) obtained from
expansion (\ref{eq:bottomexp}). Fifth and eighth columns represent
differences $F(\lambda,k)-F_1(\lambda,k)$ and
$F(\lambda,k)-F_2(\lambda,k)$, respectively.  The numbers
$\Delta_1$, $\Delta_2$ are defined in (\ref{eq:DeltaN}).}

\bigskip

Now we turn to expansion (\ref{eq:Fasymptop-final}).  From
(\ref{eq:An-represent}) we have
\[
A_0(x)=\frac{1}{\sqrt{x}}\ln(\sqrt{1+x}+\sqrt{x}),
\]
\[
A_1(x)=\frac{1}{4x}\left(\frac{1}{\sqrt{x}}\ln(\sqrt{1+x}+\sqrt{x})
-\frac{1-x}{\sqrt{1+x}}\right).
\]
Hence, the first and the second order approximations obtained from
(\ref{eq:Fasymptop-final}) are:
\begin{equation}\label{eq:F1-top}
\tilde{F}_1(\lambda,k)=K(k)-\ln\!\left(\!\!\sqrt{1+\frac{1-\lambda^2}{1-k^2}}+\sqrt{\frac{1-\lambda^2}{1-k^2}}\right),
\end{equation}
\begin{equation}\label{eq:F2-top}
\tilde{F}_2(\lambda,k)=\tilde{F}_1(\lambda,k)-\frac{1-k^2}{4}\ln\!\left(\!\!\sqrt{1+\frac{1-\lambda^2}{1-k^2}}+\sqrt{\frac{1-\lambda^2}{1-k^2}}\right)
-\frac{\lambda^2-k^2}{4\sqrt{1+\frac{1-k^2}{1-\lambda^2}}}\,.
\end{equation}
Denote by $\tilde{\Delta}_N$ the difference between the upper and
the lower bounds in (\ref{eq:error-bound1}):
\begin{multline}\label{eq:tildeDeltaN}
\tilde{\Delta}_N=\frac{(1-\lambda^2)^{N+1/2}}{2\lambda^2N\sqrt{2-k^2-\lambda^2}}
-\frac{(1-\lambda^2)^{N-1/2}(1/2)_{N}}{2NN!}
\\
\times\!\!\left(\!\! \sqrt{2-\lambda^2-k^2}
-\frac{(1-k^2)}{2\sqrt{1-\lambda^2}}
\ln\!\left\{1+2\frac{1-\lambda^2}{1-k^2}+2\frac{\sqrt{(1-\lambda^2)(2-\lambda^2-k^2)}}{1-k^2}\right\}\right).
\end{multline}
Approximation (\ref{eq:F1-top}) combined with inequality
(\ref{eq:error-bound1}) puts $F(\lambda,k)$ within an interval of
length $\tilde{\Delta}_1$, while (\ref{eq:F2-top}) puts
$F(\lambda,k)$ within an interval of length $\tilde{\Delta}_2$.
Numerical results are presented in Table~2.
{\small
$$
\begin{array}{!{\vrule width 1.2pt\relax}c|c|c!{\vrule width
1.2pt\relax}c|c|c!{\vrule width 1.2pt\relax}c|c|c!{\vrule width
1.2pt\relax}}
\hline\hline\hline\hline\hline &&&&&&&&\\[-10pt]
\lambda & k & F(\lambda,k) & \begin{array}{c}\text{1st order}\\[-3pt]\text{approx. (\ref{eq:F1-top})}\end{array} &
\begin{array}{c}\text{Absolute}\\[-3pt]\text{error}\end{array} &
\begin{array}{c}\text{Length of}\\[-3pt]\text{error range}\\[-3pt]\tilde{\Delta}_1\end{array} &
\begin{array}{c}\text{2nd order}\\[-3pt]\text{approx. (\ref{eq:F2-top})}\end{array}&
\begin{array}{c}\text{Absolute}\\[-3pt]\text{error}\end{array} &
\begin{array}{c}\text{Length of}\\[-3pt]\text{error range}\\[-3pt]\tilde{\Delta}_2\end{array}
\\[5pt]
\hline\hline\hline\hline\hline
.8 & .8   & 1.0178 & 1.1139 & -.09611 & .1509 & 1.0346 & -.01679 & .02932\\
.9 & .9   & 1.3532 & 1.3992 & -.04600 & .0576 & 1.3573 & -.00414 & .006075\\
.95 & .95 & 1.6861 & 1.7086 & -.02251 & .0252 & 1.6872 & -.00103 & .001387\\
.99 & .99 & 2.4708 & 2.4752 & -.00443 & .0045 & 2.4708 & -.408\!\!\times\!\!10^{-4} & .5164\!\!\times\!\!10^{-4}\\
.99 & .95 & 2.1496 & 2.1523 & -.00271 & .0028 & 2.1497 & -.299\!\!\times\!\!10^{-4} & .3102\!\!\times\!\!10^{-4}\\
.999 & .99& 3.0445 & 3.0447 & -.200\!\!\times\!\!10^{-3} & .200\!\!\times\!\!10^{-3}& 3.0445 & -.229\!\!\times\!\!10^{-6}&.226\!\!\times\!\!10^{-6}\\
\hline\hline\hline\hline\hline
\end{array}
$$
}\textbf{Table~2.} \textsl{Numerical examples for approximations
(\ref{eq:F1-top}) and (\ref{eq:F2-top}) obtained from expansion
(\ref{eq:Fasymptop-final}). Fifth and eighth columns represent
differences $F(\lambda,k)-\tilde{F}_1(\lambda,k)$ and
$F(\lambda,k)-\tilde{F}_2(\lambda,k)$, respectively.  The numbers
$\tilde{\Delta}_1$, $\tilde{\Delta}_2$ are defined in
(\ref{eq:tildeDeltaN}).}

\bigskip

We will compare these results with the corresponding results from
\cite{CG1} given by inequalities (\ref{eq:CG1-error}) and
(\ref{eq:CG2-error}). Denote by $\Delta_1^{\!*}$ and
$\Delta_2^{\!*}$ the interval lengthes for absolute error bounds
implied by (\ref{eq:CG1-error}) and (\ref{eq:CG2-error}),
respectively, i.e.
\begin{equation}\label{eq:delta-star}
\Delta_1^{\!*}=[\text{rhs of (\ref{eq:CG1-error})} - \text{lhs of
(\ref{eq:CG1-error})}]F(\lambda,k),~~~\Delta_2^{\!*}=[\text{rhs of
(\ref{eq:CG2-error})} - \text{lhs of
(\ref{eq:CG2-error})}]F(\lambda,k).
\end{equation}
{\small
$$
\begin{array}{!{\vrule width 1.2pt\relax}c|c|c!{\vrule width
1.2pt\relax}c|c|c!{\vrule width 1.2pt\relax}c|c|c!{\vrule width
1.2pt\relax}}
\hline\hline\hline\hline\hline &&&&&&&&\\[-10pt]
\lambda & k & F(\lambda,k) & \begin{array}{c}\text{1st order}\\[-3pt]\text{approx. (\ref{eq:CG1})}\end{array} &
\begin{array}{c}\text{Absolute}\\[-3pt]\text{error}\end{array} &
\begin{array}{c}\text{Length of}\\[-3pt]\text{error range}\\[-3pt]\Delta_1^{\!*}~\text{(see(\ref{eq:CG1-error}))}\end{array} &
\begin{array}{c}\text{2nd order}\\[-3pt]\text{approx. (\ref{eq:CG2})}\end{array}&
\begin{array}{c}\text{Absolute}\\[-3pt]\text{error}\end{array} &
\begin{array}{c}\text{Length of}\\[-3pt]\text{error range}\\[-3pt]\Delta_2^{\!*}~\text{(see(\ref{eq:CG2-error}))}\end{array}
\\[5pt]
\hline\hline\hline\hline\hline
.8 & .8   & 1.0178 & .85814 & .15968 & .2032 & .96415 & .05366 & .12508 \\
.9 & .9   & 1.3532 & 1.2278 & .12538 & .1304 & 1.3291 & .02411 & .05376 \\
.95 & .95 & 1.6861 & 1.5993 & .08687 & .0742 & 1.6771 & .00900 & .01867\\
.99 & .99 & 2.4708 & 2.4417 & .02910 & .0169 & 2.4702 & .647\!\!\times\!\!10^{-3} & .00115\\
.99 & .95 & 2.1496 & 2.0973 & .05234 & .0409 & 2.1466 & .00301 & .00898\\
.999 & .99& 3.0445 & 3.0306 & .01392 & .0076 & 3.0444 & .156\!\!\times\!\!10^{-3}&.427\!\!\times\!\!10^{-3}\\
.95 & .99 & 1.7951 & 1.7232 & .07182 & .0537 & 1.7896 & .00545 & .00750\\
.99 & .999& 2.6240 & 2.6016 & .02232 & .0115 & 2.6236 & .337\!\!\times\!\!10^{-3}& .368\!\!\times\!\!10^{-3}\\
\hline\hline\hline\hline\hline
\end{array}
$$
} \textbf{Table 3.} \textsl{Numerical examples for approximations
(\ref{eq:CG1}) and (\ref{eq:CG2}) due to Carlson and Gustafson.
Fifth and eighth columns equal $\theta_1F(\lambda,k)$,
$\theta_2F(\lambda,k)$, respectively.  The numbers
$\Delta_1^{\!*}$ $\Delta_2^{\!*}$ are defined in
(\ref{eq:delta-star}).}

\bigskip

\paragraph{6. Acknowledgements.}  The first author is partially supported by
the Russian Basic Research Fund (grant no. 05-01-00099), Far
Eastern Branch of the Russian Academy of Sciences  (grant no.
06-III-B-01-020) and INTAS (grant no.05-109-4968).


\begin{thebibliography}{99}
\bibitem{AS} M.\,Abramowitz and I.A.\,Stegun, \emph{Handbook of Mathematical
Functions}, Dover, New York, 1970.
\bibitem{Byrd} P.F.\,Byrd and M.D.\,Friedman, \emph{Handbook of elliptic
integrals for Engineers and Scienteists}, 2nd ed., Spinger-Verlag,
New York, 1971.
\bibitem{Bat1} A.\,Erd\'{e}lyi, W.\,Magnus, F.\,Oberhettinger and
F.G.\,Tricomi, \emph{Higher transcendental functions, Vol. 1},
McGraw-Hill Book Company, Inc., New York, 1953.
\bibitem{Carlson1} B.C.\,Carlson, Some series and bounds for incomplete elliptic integrals, \emph{J.
 Math. Phys.} {\bf 40}(1961), 125-134.
\bibitem{CarlsonBook} B.C.\,Carlson, \emph{Special Functions of Applied Mathematics}, Academic Press, New York, 1977.
\bibitem{CG1} B.C.\,Carlson and J.L.\,Gustafson Asymptotic
expansion of the first elliptic integral, \emph{SIAM J.
Math.Anal.}, vol.\textbf{16} (1985), no.5, 1072-1092.
\bibitem{CG2} B.C.\,Carlson and J.L.\,Gustafson,
Asymptotic approximations for symmetric elliptic integrals,
\emph{SIAM J. Math. Anal.} \textbf{25}(1994), 288-303.
\bibitem{Gustafson}J.L.\,Gustafson, \emph{Asymptotic Formulas for
Elliptic Integrals}, Ph.D. thesis, Iowa State University, Ames,
IA, 1982.
\bibitem{KSS}  D.\,Karp, A.\,Savenkova, S.M.\,Sitnik, Series expansions and
asymptotics for the third incomplete elliptic integral via partial
fraction decompositions, \emph{Journal of Computational and
Applied Mathematics}  (to appear, 2006).
\bibitem{Kaplan} E.L.\,Kaplan, Auxiliary table for the incomplete
elliptic integrals, \emph{J. Math. and Phys.} \textbf{27} (1948),
11-36.
\bibitem{Kelisky} R.P.\,Kelisky, Inverse elliptic functions and Legendre
polynomials,  \emph{Amer. Math. Monthly} 66(1959), 480-483.
\bibitem{Koepf} W.\,Koepf, \emph{Hypegeometric summation}, Advanced
Lectures in Mathematics, Vieweg, 1998.
\bibitem{Lopez} J.L.\,L\'{o}pez, Asymptotic expansions of symmetric standard
elliptic integrals, \emph{SIAM J. Math. Anal.} \textbf{31},
4(2000), 754-775.
\bibitem{Lopez1} J.L.\,L\'{o}pez, Uniform asymptotic expansions of symmetric elliptic
integrals, \emph{Constructive Approximation} \textbf{17}, 4(2001),
535-559.
\bibitem{Lopez2} J.L.\,L\'{o}pez, Asymptotic expansions of Mellin convolutions
by means of analytic continuation, \emph{J. of Comp. and Applied
Math.} (to appear, 2006).
\bibitem{Mitrinovic}D.S.\,Mitrinovic, J.E.\,Pecaric, A.M.\,Fink, \emph{Classical and new
inequalities in Analysis.} Kluwer Academic Publishers, 1993.
\bibitem{Nellis} W.J.\,Nellis and B.C.\,Carlson, Reduction and
evaluation of elliptic integrals, \emph{Math. Comp.}, \textbf{20}
(1966), 223-231.
\bibitem{Radon} B.\,Radon, Sviluppi in serie degli integrali
ellipttici, \emph{Atti. Accad. Naz. Lincei, Mem., Cl. Sci. Fis.
Mat. Nat.} Ser.(8) \textbf{2}(1950), 69-109.
\bibitem{Ponnusamy}S.\,Ponnusamy and M.\,Vuorinen, Asymptotic
expansions and inequalities for hypergeometric functions,
\emph{Mathematika} \textbf{44}(1997), 278-301.
\bibitem{Prud3} A.P.\,Prudnikov,  Yu.A.\,Brychkov and
O.I.\,Marichev, \emph{Integrals and series, Volume 3: More Special
Functions}, Gordon and Breach Science Publishers, 1990.
\bibitem{Szego} G.\,Szeg\H{o}, \emph{Orthogonal polynomials}, AMS Colloquium Publications
{\bf 23} (1991), 8th printing.
\bibitem{Sitnik1} S.M.\,Sitnik, \emph{Inequalities for the Legendre complete elliptic integrals}, Preprint of the
Far-Eastern Branch of the Russian Academy of Sciences,
Vladivostok, 1994.
\bibitem{Sitnik2} S.M.\,Sitnik, Refinements of integral Cauchy-Bunyakovskii inequality, \emph{Vestnik of Samara
State Technical University},  \textbf{9} (2000), 37-45.
\bibitem{Sitnik3} S.M.\,Sitnik, Refinements of  Cauchy-Bunyakovskii inequality and
aplications, \emph{Vestnik of Samara State Academy of Economics},
\textbf{1}, 8(2002), 302-313.
\bibitem{Sitnik4} S.M.\,Sitnik, Means and Generalizations of  Cauchy-Bunyakovskii
Inequalities and Applications, \emph{Scientific Researches of
Chernozemie} 1(2005), 3 - 42.
\end{thebibliography}
\end{document}